\documentclass[12pt]{amsart}
\usepackage[all]{xy}
\usepackage{amssymb}
\usepackage{amsfonts,amsthm,amsmath,amscd}
\begin{document}{}
\hfuzz=30pt
\newtheorem{theorem}{Theorem}[section]
\newtheorem{lemma}[theorem]{Lemma}
\newtheorem{proposition}[theorem]{Proposition}
\newtheorem{corollary}[theorem]{Corollary}

\theoremstyle{definition}
\newtheorem{definition}[theorem]{Definition}
\newtheorem{example}[theorem]{Example}
\newtheorem{xca}[theorem]{Exercise}

\newtheorem{remark}[theorem]{Remark}

\numberwithin{equation}{section}

\newcommand{\V}{{\mathcal{V}=(\mathcal{V}_{0}, \otimes, I)}}
\newcommand{\VV}{{\mathcal {V}_0}}
\newcommand{\M}{{\text{Mod}}}

\title{On comonadicity of the extension-of-scalars functors}


\author{Bachuki Mesablishvili}
\address{Razmadze Mathematical Institute, Tbilisi 0193, Republic of Georgia}
\email{bachi@rmi.acnet.ge}




\begin{abstract} A criterion for comonadicity of the extension-of-
scalars functor associated to an extension of (not necessarily
commutative) rings is given. As an application of this criterion,
some known results on the comonadicity of such functors are
obtained.

\end{abstract}

\maketitle

\section{Introduction}

In view of the observation of Caenepeel (see \cite{C}) that
noncommutative descent for modules reduces to comonadicity of the
corresponding extension-of-scalars functor, it becomes even more
sensible to have manageable tests for comonadicity of the
extension-of-scalars functors. (Although there are several results
obtained along these lines (see \cite{Br}, \cite{Br1}, \cite{C},
\cite{C1}, \cite{Ci}, \cite{EG}, \cite{JT}, \cite{M}, \cite{Nu})
the question of comonadicity of such functors is not fully
answered yet.) The main result of this note gives such a test.

For the basic definitions of category theory, see \cite{Mc}

\section{Preliminaries}

A monoidal category $\V$ is called \emph{biclosed} if for all $X
\in \text{Ob}(\VV)$, the functors $$-\otimes X, \,\,X \otimes - :
\VV \to \VV$$ have (chosen) right adjoints, denoted $[X, -]$ and
$\{X, -\}$, respectively. In other words, a biclosed monoidal
category consists of a monoidal category $\V$, equipped with two
functors $$[-, -]\, \, \, ,\{-, -\} : \VV ^{op}\times \VV \to
\VV,$$ for which there are natural isomorphisms
\begin{equation}
\VV (X, [Y, Z]) \simeq \VV(X\otimes Y, Z)\simeq \VV (Y, \{X, Z\})
\end{equation}

Recall that the adjunctions $X \otimes - \dashv \{X, -\}$ and $-
\otimes X \dashv [X, -]$ are internal, in the sense that one has
natural isomorphisms
\begin{equation}
\{X \otimes Y, Z\} \simeq \{Y, \{X, Z\}\}
\end{equation}
and
\begin{equation}
[X \otimes Y, Z] \simeq [X, [Y, Z]].
\end{equation}

\smallskip
Let us recall that a morphism in a category $\mathcal{A}$ is a
\emph{regular monomorphism} if it is an equalizer of some pair of
morphisms. Recall also that an object $X$ of $\mathcal{A}$ is
\emph{injective} if it is injective with respect to the class of
regular monomorphisms of $\mathcal{A}$, that is, if every
extension problem

$$\xymatrix{
A \ar[d]_f \ar @{^{(}->}[r]^m &B \ar@{.>}[dl]^{\bar{f}}\\
X}
$$ with $m$ a regular monomorphism has a solution $\bar{f}: B \to
X$ extending $f$ along $m$, i.e., satisfying $\bar{f}m=f$. (The
dual notions are the \emph{regular epimorphism} and the
\emph{projective object}.)

\smallskip

Let $\V$ be a monoidal category and let $f : X \to Y$ be a
morphism in $\VV$. We say that $f$ is \emph{right} (resp.
\emph{left}) \emph{pure} if, for any $Z \in \text{Ob}(\VV)$, the
morphism $$f \otimes Z : X \otimes Z \to Y \otimes Z)$$
$$(\text{resp.} \,\, Z \otimes f : Z\otimes X \to Z \otimes Y$$
is a regular monomorphism.

\bigskip

Henceforth, we suppose without explicit mention that $\mathcal{V}$
is a finitely complete and finitely cocomplete  monoidal biclosed
category whose unit $I$ for the tensor product is projective.

\smallskip

\begin{theorem} Let $Q$ be an object of $\VV$ for which the functor
$$\{-, Q\}: \VV ^{op} \to \VV$$ is conservative
(that is, isomorphism-reflecting) and preserves regular
epimorphisms. Then the following properties of a morphism $f : X
\to Y$ of \, $\VV$ are equivalent:
\begin{itemize}
\item[(i)] The morphism $f$ is right pure.

\item[(ii)] The morphism $f \otimes \{X, Q\} : X \otimes \{X, Q\}
\to Y \otimes \{X, Q\}$ is a regular monomorphism.

\item[(iii)] The morphism $\{f, Q\} : \{Y, Q\} \to \{X, Q\}$ is a
split epimorphism.
\end{itemize}
\end{theorem}
\begin{proof} $\text{(i)}$ implies $\text{(ii)}$ trivially. To see that
$\text{(ii)}$ implies $\text{(iii)}$, let us assume that the
morphism
$$f \otimes \{X, Q\}: X \otimes \{X, Q\} \to Y \otimes \{X, Q\}$$
is a regular monomorphism. Since the functor $\{-, Q\}: \VV ^{op}
\to \VV$ preserves regular epimorphisms by hypothesis, the
morphism
$$\{\{X, Q\},\{f, Q\}\}: \{\{X, Q\}, \{Y, Q\}\} \to \{\{X, Q\},
\{X, Q\}\},
$$ which is isomorphic by $(2.2)$ to the morphism $$\{f \otimes
\{X, Q\}, Q\} : \{Y \otimes \{X, Q\}, Q\} \to \{X \otimes \{X,
Q\}, Q\},$$ is a regular epimorphism in $\VV$. Since $I$ is
assumed to be projective in $\VV$, the functor $$\VV(I, -) : \VV
\to \bold{Set}$$ takes regular epimorphisms to surjections. It
follows that the map $$\VV( \{X, Q\},\{f, Q\})$$ of sets, which
(using (2.1))) is isomorphic to the map
$$ \VV (I, \{\{X, Q\} ,\{f, Q\}\}),$$  is surjective. But this
means that every morphism
$$\{X, Q\} \to \{X, Q\}$$ factors through $\{f, Q\}$, that is to
say, that $\{f, Q\}$ is a split epimorphism, as is seen from the
special case of the identity morphism $1_{\{X,\,\, Q\}}$.

It remains to show that $\text{(iii)}$ implies $\text{(i)}$. If
the morphism $$\{f, Q\} : \{Y, Q\} \to \{X, Q\}$$ is a split
epimorphism, then so is
$$\{Z,\{f, Q\}\}: \{Z, \{Y, Q\} \} \to \{Z, \{X, Q\}\}$$ too, for all
$Z \in \text{Ob} (\VV)$. Identifying the morphism $\{Z,\{f, Q\}\}$
(via the isomorphism $(2.2)$) with $\{f \otimes Z, Q\}$, we see
that the morphism
$$\{f \otimes Z, Q\} : \{Y \otimes Z ,Q\} \to \{X \otimes Z, Q\}$$
is also a split epimorphism. We now observe that, since the
functor $\{-, Q\}: \VV ^{op} \to \VV$ admits as a right adjoint
the functor $[-, Q]: \VV  \to \VV ^{op}$, as can be seen from the
following sequence of natural isomorphisms: $$\VV (X, \{Y,
Q\})\simeq \VV (Y\otimes X, Q) \simeq \VV (Y, [X, Q])\simeq \VV
^{op}([X, Q], Y),$$ to say that $\{-, Q\}$ is conservative and
preserves regular epimorphisms is to say that it preserves and
reflects regular epimorphisms. And since any split epimorphism is
regular, it follows that the morphism $f \otimes Z : X \otimes Z
\to Y \otimes Z$ is a regular monomorphism for all $Z \in
\text{Ob} (\VV)$. Thus $\text{(iii)}$ implies $\text{(i)}$. The
proof of the theorem is now complete.
\end{proof}

There is of course a dual result:

\begin{theorem} Let $Q$ be an object of $\VV$ such that the functor
$$[-, Q]: \VV ^{op} \to \VV$$ is conservative and preserves regular
epimorphisms. Then the following properties of a morphism $f : X
\to Y$ of \,$\VV$ are equivalent:
\begin{itemize}
\item[(i)] The morphism $f$ is left pure.

\item[(ii)] The morphism $ [X, Q] \otimes f : [X, Q] \otimes X \to
[X, Q] \otimes Y$ is a regular monomorphism.

\item[(iii)] The morphism $[f, Q ]: [Y, Q] \to [X, Q]$ is a split
epimorphism.
\end{itemize}
\end{theorem}

An object $Q$ of a monoidal biclosed category
$$\mathcal{V}=(\mathcal{V}_{0}, \otimes, I, [-,-], \{-,-\})$$ is
said to be \emph{cyclic} if the functors $\{-, Q\}$ and $[-, Q]$
are naturally isomorphic. If $Q$ is such an object, we shall
denote by $[[-, Q]]$ the functor $[-,Q]\simeq \{-,Q\}$.

\smallskip
\bigskip

Combining Theorems 2.1 and 2.2, we get:

\begin{theorem}Let $Q$ be a cyclic object of $\VV$ for which the functor
$$[[-, Q]]: \VV ^{op} \to \VV$$ is conservative and preserves regular
epimorphisms (equivalently, preserves and reflects regular
epimorphisms). Then the following properties of a morphism $f : X
\to Y$ of \,$\VV$ are equivalent:
\begin{itemize}
\item[(i)] The morphism $f$ is left pure.

\item[(ii)] The morphism $f$ is right pure.

\item[(iii)] The morphism $ [[X, Q]] \otimes f : [[X, Q]] \otimes
X \to [[X, Q]] \otimes Y$ is a regular monomorphism.

\item[(iv)] The morphism $ f \otimes [[X, Q]] : X \otimes [[X, Q]]
\to Y \otimes [[X, Q]]$ is a regular monomorphism.

\item[(v)] The morphism $[[f, Q]] : [[Y, Q]] \to [[X, Q]]$ is a
split epimorphism.
\end{itemize}
\end{theorem}

\bigskip

\section{A Criterion for Comonadicity of Extension-of-Scalars Functors}

In this section we present our main result.

Let us fix a commutative ring $K$ with unit ($K=\mathbb{Z}$, the
ring of integers, inclusive). All rings under consideration are
associative unital $K$-algebras. A right or left module means a
unital module. All bimodules are assumed to be $K$-symmetric. The
$K$-categories of left and right modules over a ring $A$ are
denoted by ${_A \M}$ and $\M _A$, respectively; while the category
of $(A,B)$-bimodules is $_A \M _B$. We will use the notation ${_B
M} _A$ to indicate that $M$ is a left $B$, right $A$-module.

It is a well-known fact that, for a fixed ring $A$, the category
$_A \M _A$ is a monoidal category with tensor product of two
$(A,A)$-bimodules being their usual tensor product over $A$ and
the unit for this tensor product being the $(A,A)$-bimodule $A$.
Moreover, this monoidal category is biclosed: If $M$ and $N$ are
two $(A,A)$-bimodules, then $[M, N]=\M _A (M, N)$ and $\{M,
N\}={_A \M} (M, N)$.

\smallskip
\smallskip

For any $(A,A)$-bimodule $M$, the \emph{character}
$(A,A)$-\emph{bimodule} of $M$ is defined to be $M^+
=\textbf{Ab}(M, \mathbb{Q} /\mathbb{Z})$ (where $\textbf{Ab} $ is
the category of abelian groups and $\mathbb{Q} /\mathbb{Z}$ is the
rational circle abelian group). This is an $(A,A)$-bimodule via
the actions $(afa')(m)=f(a'ma)$.

\begin{lemma}The character bimodule $A^+$ of the $(A,A)$-bimodule
$A$ is a cyclic object of the monoidal biclosed category ${_A
\emph{\M}} _A$ of $(A,A)$-bimodules.
\end{lemma}

\begin{proof} The following string of natural isomorphisms
$$\{-, A^+\}=\{-,\textbf{Ab}(A, \mathbb{Q} /\mathbb{Z})\}=
{_A \M (-, \textbf{Ab}(A, \mathbb{Q} /\mathbb{Z}))}\simeq$$
$$\simeq \textbf{Ab}(A \otimes_A -, \mathbb{Q} /\mathbb{Z})
\simeq \textbf{Ab}( -, \mathbb{Q} /\mathbb{Z})\simeq \textbf{Ab}(-
\otimes_A A, \mathbb{Q} /\mathbb{Z})\simeq $$$$\simeq \M _A (-,
\textbf{Ab}(A, \mathbb{Q} /\mathbb{Z}))=[-, A^+]$$ shows that the
functors $$\{-, A^+\}, \,\, [-, A^+]: ({_A \M} _A)^{op}\to {_A \M}
_A$$ are naturally equivalent.
\end{proof}

Since the functor $[[-, A^+]]$ is naturally equivalent to
$\textbf{Ab}(-, \mathbb{Q} /\mathbb{Z})$ and since $\mathbb{Q}
/\mathbb{Z}$ is an injective cogenerator in $\textbf{Ab}$, we have
that

\begin{lemma}The functor $$[[-, A^+]]: ({_A \emph{\M}} _A)^{op}\to {_A
\emph{\M}} _A$$ is exact and conservative.
\end{lemma}

Before we prove our main result we recall a result from \cite{M}:

\begin{theorem} Let $i: A \to B$ be a homomorphism of rings. If the
induced morphism $i^+ : B^+ \to A^+$ is a split epimorphism of
$(A,A)$-bimodules, then the functors $$- {\otimes_A B} :
\emph{\M}_A \to \emph{\M}_B$$ and $${B \otimes_A} - : {_A
\emph{\M}} \to {_B \emph{\M}}$$ are both comonadic.
\end{theorem}

Recall (for example from \cite{L}) that a morphism $ f : M \to N$
of right $A$-modules is called \emph{pure} if $ f \otimes _A 1_L :
M \otimes_A L \to N \otimes _A L$ is injective for every left
$A$-module $L$. Pure morphisms in the category of left $A$-modules
are defined analogously.

\smallskip
\smallskip
The main result of this note is contained in the following:

\begin{theorem} Let $i: A \to B$ be a homomorphism of rings. If $A$
is a separable $K$-algebra, then the following are equivalent:

\begin{itemize}
\item[(i)] $i$ is a pure morphism of left $A$-modules.

\item[(ii)] $i$ is a pure morphism of right $A$-modules.

\item[(iii)] $i^+ : B^+ \to A^+$ is a split epimorphism of
$(A,A)$-bimodules.

\item[(iv)] The functor $- {\otimes_A B} : \emph{\M}_A \to
\emph{\M}_B$ is comonadic.

\item[(v)] The functor ${B \otimes_A} - : {_A \emph{\M}} \to {_B
\emph{\M}}$ is comonadic.
\end{itemize}
\end{theorem}

\begin{proof}We remark first that, by
left-right symmetry, it suffices to prove the equivalence of
$\text{(i)}$, $\text{(iii)}$ and $\text{(v)}$.

Using that the forgetful functor ${_A \M _A} \to \M_A$ preserves
and reflects monomorphisms and tensor products, it is easy to see
that if $i$ is a pure morphism of left $A$-modules, then it is
left pure in ${_A \M _A}$. And since assuming $A$ be $K$-separable
is, just by definition, the same as assuming $A$ be projective in
${_A \M_A }$, it follows from Theorem 2.3 that the morphism $[[i,
A^+]] \simeq i^+$ is a split epimorphism of $(A,A)$-bimodules,
provided that $i$ is a pure morphism of left $A$-modules. Thus
$\text{(i)}$ implies $\text{(iii)}$.

$\text{(iii)}$ implies $\text{(v)}$ by Theorem 3.3.

It is well known that the functor $${-\otimes_A B} : \M_A \to
\M_B$$ admits as a right adjoint the functor $$\M_B (B, -) : \M_B
\to \M_A$$ and that the unit $\eta$ of this adjunction  has
components $$\eta_X : X \otimes_A i : X \simeq X \otimes_A A  \to
X \otimes_A B, \,\, X \in \M_A .$$ Thus $i$ is a pure morphism of
left $A$-modules precisely when $\eta$ is componentwise a
monomorphism. According to Theorem 9 of Section 2.3 of \cite{BW},
this is in particular the case when the functor ${-\otimes_A B}$
is comonadic. So $\text{(v)}$ implies $\text{(i)}$. This completes
the proof of the theorem.
\end{proof}

\section{Applications}

In this section we state some consequences of our main theorem. To
state the first one, we need a definition. Let $A$, $B$ be rings.
Recall \cite{C1} that an $(A,B)$-bimodule $M$ is said to be
\emph{totally faithful} as a left $A$-module if the morphism $$X
\to \M_B (M, X \otimes_A M), \,\,\, m \to x \otimes_A m,$$ is
injective for every $X \in \M_A$, or equivalently, if the unit of
the adjunction
$${-\otimes_A M} \dashv \M_B (M, -) : \M_B \to \M_A$$ is pointwise
a monomorphism.

\begin{theorem}Let $A$ and $B$ be rings, $M$ an $(A,B)$-bimodule
with $M_B$ finitely generated and projective, $\mathcal E _M
=\emph{\M}_B (M, M)$ the right endomorphism ring of $M_B$ and
$$i_M : A \to \mathcal E _M, \,\, a \to [m \to am]$$ the corresponding
ring homomorphism. If $A$ is $K$-separable, then the following are
equivalent:
\begin{itemize}

\item[(i)] The bimodule ${_A M_B}$ is totally faithful as a left
$A$-module.

\item[(ii)] The bimodule ${_B {M^*}_A}$ is totally faithful as a
right $A$-module. (Here we denote by $M^*$ the dual $\emph{\M}_B
(M, B)$ of $M_B$ which is a $(B,A)$-bimodule in a canonical way.)

\item[(iii)] The morphism $(i_M)^+ : (\mathcal E _M)^+ \to A^+$ is
a split epimorphism of $(A,A)$-bimodules.

\item[(iv)] The functor ${-\otimes_A M}: \emph{\M}_A \to
\emph{\M}_B$ is comonadic.

\item[(iv)] The functor ${{M^*} \otimes_A -}: {_A \emph{\M}} \to
{_B \emph{\M}}$ is comonadic.

\end{itemize}

\end{theorem}

\begin{proof} Immediate from Theorem 3.4 using that:

\begin{itemize}
\item ${_A M_B}$ (resp. ${_B M^* _A}$) is totally faithful as a
left (resp. a right) $A$-module if and only if $i_M : A \to
\mathcal E _M$ is a pure morphism of left (resp. right)
$A$-modules (see Lemma 2.2 in \cite{C1}, or Proposition 7.3 in
\cite{M});

\item the functor ${-\otimes_A M}: \M_A \to \M_B$ (resp. ${{M^*}
\otimes_A -}: {_A\M} \to {_B\M}$) is comonadic if and only if the
functor ${-\otimes_A \mathcal E _M}: \M_A \to \M_{\mathcal E _M}$
(resp. ${{\mathcal E _M} \otimes_A -}: {_A\M} \to {_{\mathcal E
_M}\M}$) is so (see Theorem 7.5 in \cite{M}).
\end{itemize}
\end{proof}

As a special case of Theorem 4.1 one can take $A=K$. Then, since
obviously $K$ is $K$-separable, we recover a result by Caenepeel,
De Groot and Vercruysse \cite{C1}.

\begin{theorem} [Caenepeel, De Groot
and Vercruysse \cite{C1}]Let $A$ be a ring and let $M$ be a
$(K,A)$-bimodule with $M_A$ finitely generated and projective.
Then the following are equivalent:
\begin{itemize}

\item[(i)] The morphism $i_M : K \to \mathcal E _M = \emph{\M}_A
(M, M)$ is a pure morphism of left $K$-modules.

\item[(ii)] The morphism $i_M : K \to \mathcal E _M = \emph{\M}_A
(M, M)$ is a pure morphism of right $K$-modules.

\item[(iii)] The bimodule ${_A M_B}$ is totally faithful as a left
$A$-module.

\item[(iv)] The bimodule ${_B {M^*}_A}$ is totally faithful as a
right $A$-module.

\item[(v)] The functor ${-\otimes_A M}: \emph{\M}_A \to
\emph{\M}_B$ is comonadic.

\item[(vi)] The functor ${{M^*} \otimes_A -}: {_A \emph{\M}} \to
{_B \emph{\M}}$ is comonadic.

\end{itemize}
\end{theorem}

For the special case in which $M=A$, we recapture easily the
following result of  Joyal and Tierney (unpublished, but see
\cite{M1}). Recall (for example from\cite{JT}) that a homomorphism
$i: K \to A$ of commutative rings is said to be \emph{effective
for descent} if the extension-of-scalars functor
$$ A \otimes_K - : \M_K \to \M_A$$ is comonadic.

\begin{theorem} [Joyal and Tierney]A homomorphism $i: K \to A$
of commutative rings is effective for descent if and only if it is
a pure morphism of (say left) $K$-modules.
\end{theorem}

We end this note with an interesting consequence of Theorem 3.4.
Let us write $\mathbb{M}_n(K)$ for the ring of $n\times n$
matrices over $K$.

\begin{theorem} The following are equivalent for a homomorphism
$i: \mathbb{M}_n(K) \to A$ of rings:

\begin{itemize}

\item[(i)] $i$ is a pure morphism of left $A$-modules;

\item[(ii)] $i$ is a pure morphism of right $A$-modules;

\item[(iii)] the functor $- {\otimes_{\mathbb{M}_n(K)} A} :
\emph{\M}_{\mathbb{M}_n(K)} \to \emph{\M}_A$ is comonadic;

\item[(iv)] the functor ${A \otimes_{\mathbb{M}_n(K)}} - :
{_{\mathbb{M}_n(K)} \emph{\M}} \to {_A \emph{\M}}$ is comonadic.
\end{itemize}
\end{theorem}

\begin{proof}Immediate from Theorem 3.4, since for any
$n \in \mathbb{N}$, the ring $\mathbb{M}_n(K)$ is $K$-separable.
\end{proof}

\bigskip

\bibliographystyle{amsplain}

\end{document}